\documentclass[10pt]{amsart}
\usepackage{amssymb}
\usepackage{epsfig}
\usepackage{url}
\usepackage{setspace}
\theoremstyle{plain}

\newtheorem{thm}{Theorem}[section]

\newtheorem{rem}[thm]{Remark}
\newtheorem{ques}[thm]{Question}

\def\cal{\mathcal}
\def\bbb{\mathbb}
\def\op{\operatorname}
\renewcommand{\phi}{\varphi}

\newcommand{\N}{\bbb{N}}
\newcommand{\Z}{\bbb{Z}}

\begin{document}

\title[On certain diophantine equations]{On certain diophantine equations related to triangular and tetrahedral numbers}
\author{Maciej Ulas}

\bigskip

\begin{abstract}
In this paper we give solutions of certain diophantine equations
related to triangular and tetrahedral numbers and propose several
problems connected with these numbers.

The material of this paper was presented in part at the 11th
International Workshop for Young Mathematicians - NUMBER THEORY,
Krak\'{o}w, 14th-20th september 2008.
\end{abstract}

\maketitle

\section{Introduction}\label{Section1}
By a {\it triangular number} we call the number of the form
\begin{equation*}
t_{n}=1+2+\ldots + n-1 + n=\frac{n(n+1)}{2},
\end{equation*}
where $n$ is a natural number. The number $t_{n}$ can be interpreted
as a number of circles necessary to build an equilateral triangle
with side of length $n$. In analogous manner we define tetrahedral
number $T_{n}$, which gives number of balls necessary to to build
tetrahedron with side of length $n$. More explicite value of $T_{n}$
is given by
\begin{equation*}
T_{n}=t_{1}+t_{2}+\ldots+t_{n-1}+t_{n}=\frac{n(n+1)(n+2)}{6}.
\end{equation*}

Wac{\l}aw Sierpi\'{n}ski in the booklet \cite{Sier} and in the
papers \cite{Sier1,Sier2,Sier3,Sier4,Sier5} gave many interesting
results concerning the problem of solvability of diophantine
equations related to triangular nad tetrahedral numbers. The aim of
this paper is to give solutions of certain diophantine problems
which was left as open in the booklet \cite{Sier} and give some new
results. Our proofs are of elementary character and we do not assume
any special knowledge from number theory.

The number of possible problems which can be stated in connection
with triangular and tetrahedral numbers is bounded only by
imagination and with interest of researcher. This is the reason of
selection of problems in this paper. We encourage the reader to
solve problems mentioned in this paper and to state own problems.

\section{Triangular and tetrahedral numbers as sides of right triangles}\label{Section2}
In this section we are interested in the construction of right
triangles with such a property that  are triangular numbers.

We start with the problem related to the construction of right
triangles with legs which are triangular numbers. So, we will be interested in integers solutions of the diophantine equation
\begin{equation}\label{protr}
z^2=t_{x}^2+t_{y}^2,
\end{equation}
W. Sierpi\'{n}ski in \cite[page 34]{Sier} has shown that the above
equation has infinitely many solutions in integers. However, all
solutions presented by him satisfied the condition
$\op{GCD}(t_{x},t_{y})>1$. In other words, triple
$X=t_{x},\;Y=t_{y},\;Z=z$ which satisfy the equation $X^2+Y^2=Z^2$
is not primitive solution. A. Schinzel showed that bthe set of
integer solutions of the equation (\ref{protr}) which satisfy the
condition $\op{GCD}(t_{x},t_{y})=1$ is infinite. It is natural to
ask the following question: does the equation (\ref{protr}) have
parametric solutions? In other words: does the equation
(\ref{protr}) have solutions in the ring $\Z[u]$?

Generally, these sort of questions is very difficult and we do not
have any general theory which could by used. However, as we will
see, for our particular equation it is possible to construct
infinitely many polynomials $x(u),y(u), z(u)\in\Z[u]$ which satisfy
the condition $\op{GCD}(t_{x(u)},t_{y(u)})=1$.

It is clear that we can consider the equation
\begin{equation*}
r^2=(p^2-1)^2+(q^2-1)^2.
\end{equation*}
Indeed, if $p,q,r$ satisfied the above equation, then the triple of integers
$(p-1)/2,(q-1)/2,r/8$ will be solution of (\ref{protr}). From this remark we can see that the quantities  $p^2-1,q^2-1,r$ must be solutions of the equation of Pythagoras $Z^2=X^2+Y^2$. It is well known that, all solutions of this equation are of the form
\begin{equation*}
X=2abc,\quad Y=(a^2-b^2)c,\quad Z=(a^2+b^2)c,
\end{equation*}
where $a,b,c$ are certain integers. Let us put $c=1$ and consider the system of equations given by
\begin{equation*}
p^2=2ab+1,\quad q^2=a^2-b^2+1.
\end{equation*}
First equation of the above system will be satisfied if we put
\begin{equation*}
a=\frac{u(ku-2)}{2},\quad b=k,\quad p=ku-1.
\end{equation*}
We put the quantities given above into the equation $q^2=a^2-b^2+1$ and we get
\begin{equation*}
q^2=\frac{u^4-4}{4}k^2-u^3k+u^2+1=:f(k).
\end{equation*}
This is Pell type equation depending on the parameter $u$. Let us
note that $f(1)=(u(u-2)/2)^2$ and that the following identity holds
\begin{equation*}
f\Big(\frac{(u^4-2)k}{2}+u^2q-u^3\Big)-\Big(\frac{u^2(u^4-4)k+2(u^4-2)q-2u^5}{4}\Big)^2=f(k)-q^2.
\end{equation*}
From the above we can deduce that if
we define
\begin{equation}\label{sys}
\begin{cases}
&k_{0}=1,\quad q_{0}=u(u-2)/2,\\
&k_{n}=\frac{(u^4-2)k_{n-1}}{2}+u^2q_{n-1}-u^3,\\
&q_{n}=\frac{u^2(u^4-4)k_{n-1}+2(u^4-2)q_{n-1}-2u^5}{4},
\end{cases}
\end{equation}
then, the polynomials
$p_{n}(u)=k_{n}(u)u-1,\;q_{n}(u),\;r_{n}(u)=u^2(k_{n}(u)-2)^2/4+k_{n}(u)^2$
for $n=1,2,\ldots$ satisfy the equation $(p^2-1)^2+(q^2-1)^2=r^2$.
Finally, we get that the polynomials
$x_{n}(u)=(p_{n}(u)-1)/2,\;y_{n}(u)=(q_{n}(u)-1)/2,\;z_{n}(u)=r_{n}(u)/8$
satisfied the equation $t_{x}^2+t_{y}^2=z^2$. In particular, for
$n=1$ we get
\begin{align*}
&x_{1}(u)=(u^5-2u^4-u-2)/2,\\
&y_{1}(u)=(u^2+1)(u^4-2u^3-u^2+2u-2)/4,\\
&z_{1}(u)=(u^{12}-4u^{11}+4u^{10}+2u^8-16u^7+24u^6-7u^4+20u^3+4u^2+4)/32.
\end{align*}
Resultant $\op{Res}(t_{x_{1}},t_{y_{1}})$ of the polynomials
$t_{x_{1}(u)},t_{y_{1}(u)}$ is equal to $2^{-58}$, which means that
the polynomials are co-prime. Let us note that the polynomials
$x_{n}(2u+1),y_{n}(2u+1),z_{n}(2u+1)$ belong to $\Z[u]$. It is
possible to prove (we will not do this here) that for each positive
integer $n$ thew polynomials $t_{x_{n}(u)},\;t_{y_{n}(u)}$ are
co-prime. We have proved the following
\begin{thm}
The equation $z^2=t_{x}^2+t_{y}^2$ has infinitely many
solutions in polynomials $x(u),y(u),z(u)\in \Z[u]$.
\end{thm}

\begin{rem}
{\rm The problem of construction of the right angle triangles which all sides are triangular numbers so the problem of the construction of integer solutions of the equation $t_{x}^2+t_{y}^2=t_{z}^2$  was posed by
K. Zarankiewicz. Only one integer solution of this equation is known so far: $x=132,y=143,z=164$.
In connection with this problem we prove the following
\begin{thm}
The equation $t_{x}^2+t_{y}^2=t_{z}^2$ has infinitely many
solutions in rational numbers.
\end{thm}
\begin{proof}
Because we are interested in rational solution of our equation, so
without loss of generality we can consider the equation
\begin{equation*}
(p^2-1)^2+(q^2-1)^2=(r^2-1)^2.
\end{equation*}
Let $u,v$ be indeterminate parameters. Let us put
\begin{equation*}
p=2uvT-1,\quad q=(v^2-u^2)T+1,\quad r=(v^2+u^2)T+1.
\end{equation*}
For the quantities $p,q,r$ defined above we have
\begin{align*}
(p^2-1)^2&+(q^2-1)^2-(r^2-1)^2=\\
         &-8T^3u^2(u+v)^2(u^2-2uv+3v^2)-8T^4
         u^2(u+v)^2(u^2v^2-2uv^3+v^4).
\end{align*}
The polynomial on the right side of the above equality has two rational roots: $T=0$ and
\begin{equation*}
T=\frac{-u^2+2uv-3v^2}{(u - v)^2v^2}.
\end{equation*}
Using now the quantity $T$, definition of $p,q,r$ and remembering
that $t_{-x}=t_{x-1}$ we find rational parametric solution of the
equation $t_{x}^2+t_{y}^2=t_{z}^2$ in the form
\begin{align*}
&x(u,v)=\frac{u(u^2-2uv+3v^2)}{(u - v)^2v},\\
&y(u,v)=\frac{(u+v)(u^2-2uv+3v^2)}{2(u - v)v^2},\\
&z(u,v)=\frac{u^4-2u^3v+2u^2v^2+2uv^3+v^4}{2(u - v)^2v^2}.
\end{align*}

Let us note that the solution we have obtained can be used in order
to prove that our equation has infinitely many solutions in
$\cal{S}$-integers, where $\cal{S}$ is finite set of primes and
$2\in\cal{S}$. Let us recall that we say that rational number $r$ is
$\cal{S}$-integer, if the set of prime divisors of the denominator
of $r$ is contained in  $\cal{S}$. Indeed, if we put
$v=V^{n},u=U^{m}-V^{n}$, where $U,V$ are finite products of the
elements from $\cal{S}$ and $2\in\cal{S}$ we obtain rational numbers
$x(u,v),y(u,v),z(u,v)$ which are $\cal{S}$-integers. The quantity of
$\cal{S}$-integers solutions of our equation may suggest that the
set of the integer solutions of the equation
$t_{x}^2+t_{y}^2=t_{z}^2$ should be bigger (then 1).
\end{proof}
}\end{rem}

Now we will consider the problem of construction of right angle triangles which legs are tetrahedral numbers. So, we are interested in the integer solutions of the diophantine equation
\begin{equation}\label{protet}
z^2=T_{x}^2+T_{y}^2.
\end{equation}
Let us note that $91^2=T_{5}^2+T_{7}^2$, so this equation has
integer solution. W. Sierpi\'{n}ski in \cite[page 57]{Sier} wrote
that it is unclear if the equation (\ref{protet}) has infinitely
many solutions in integers. However, without much trouble we can
construct infinitely many solutions ot this equation which satisfy
the condition $y-x=1$. Indeed, we have
$T_{6x}^2+T_{6x+1}^2=(3x+1)^2(6x+1)^2f(x)$, where $f(x)=8x^2+4x+1$.
Because $f(0)=1$ and the following identity holds
\begin{equation*}
f(17x + 6z+4)-(48x+17z+12)^2=f(x)-z^2,
\end{equation*}
we can see that the equation $f(x)=z^2$ has infinitely many solutions
$x_{n},z_{n}$ given by
\begin{equation*}
x_{0}=0,\quad z_{0}=1,\quad x_{n}=17x_{n-1}+6z_{n-1}+4,\quad
z_{n}=48x_{n-1}+17z_{n-1}+12.
\end{equation*}
From the above we can see that for each $n$ we get the identity
\begin{equation*}
((3x_{n}+1)(6x_{n}+1)z_{n})^2=T_{6x_{n}}^2+T_{6x_{n}+1}^2.
\end{equation*}
In particular $T_{60}^2+T_{61}^2=54839^2,\; T_{2088}^2+T_{2089}^2=
2150259925^2,\;\ldots\;.$

In the light of the above result it is interesting to ask the
question if it is possible to find infinite family of solutions
$x_{n},y_{n},z_{n}$ of the equation (\ref{protet}), with such a
property that $y_{n}-x_{n}\rightarrow\infty$?

We will construct two families which satisfied mentioned condtion.

For the proof let us put
\begin{align*}
&x(u,v)=v^2-u^2-1,\\
&y(u,v)=\frac{3v^2-2uv+3u^2-3}{2},\\
&z(u,v)=\frac{(v^2-u^2)Z(u,v)}{192},
\end{align*}
where
\begin{equation*}
Z(u,v)=105v^4-108uv^3+(150u^2-96)v^2-4u(27u^2-16)v+3(u^2-1)(35u^2+3).
\end{equation*}

For such  $x,y,z$ we get the following equality
\begin{equation*}
T_{x}^2+T_{y}^2-z^2=\frac{h(u,v)(h(u,v)+2)H(u,v)}{36864},
\end{equation*}
where $h(u,v)=-1+u^2-6uv+v^2$ and $H(u,v)$ is the polynomial of degree
8. Let us note that the equation $h(u,v)=0$ has infinitely many solutions in positive integers. This is an immediate consequence of the equality $h(6,35)=0$ and the identity
\begin{equation*}
h(u,v)=h(v,6v-u).
\end{equation*}
This means that for the sequence defined recursively by the
equations
\begin{equation*}
u_{0}=6,\quad u_{1}=35,\quad u_{n}=6u_{n-1}-u_{n-2},
\end{equation*}
we have equality $h(x_{n-1},x_{n})=0$ for $n=1,2,\ldots\;$. We
conclude that the numbers
$x(u_{n-1},u_{n}),y(u_{n-1},u_{n}),z(u_{n-1},u_{n})$ are integer
solutions of the diophantine equation $T_{x}^2+T_{y}^2=z^2$. In
particular
\begin{equation*}
T_{1188}^2+T_{1680}^2=839790700^2,\quad
T_{40390}^2+T_{57120}^2=32946833683400^2,\ldots,\;.
\end{equation*}

Let us note that in order to prove above result we can also take
\begin{align*}
&x'(u,v)=x(u,v),\\
&y'(u,v)=y(u,v)+1,\\
&z'(u,v)=\frac{(v^2-u^2)Z'(u,v)}{192},
\end{align*}
where
\begin{equation*}
Z'(u,v)=105v^4-108uv^3+(150u^2+96)v^2-4u(27u^2+16)v+3(u^2+1)(35u^2-3).
\end{equation*}
For $x',y',z'$ defined in this way we have an identity
\begin{equation*}
T_{x'}^2+T_{y'}^2-z'^2=\frac{h(u,v)(h(u,v)+2)H'(u,v)}{36864},
\end{equation*}
where $h(u,v)$ is the same polynomial we have obtained previously
and $H'$ is a certain polynomial of degree 8.

We have proved
\begin{thm}
\begin{enumerate}
\item The equation $T_{x}^2+T_{y}^2=z^2$ has infinitely many
integer solutions satisfied the condition $y-x=1$.

\item There exists an infinite sequence $(x_{n},y_{n},z_{n})$ of
solutions of the equation $T_{x}^2+T_{y}^2=z^2$ with such a property
that $y_{n}-x_{n}\rightarrow \infty$.
\end{enumerate}
\end{thm}

It is easy to see that the solutions of the equation (\ref{protet})
we have obtained are not co-prime. This suggest the following:

\begin{ques}
Does the equation $z^2=T_{x}^2+T_{y}^2$ have infinitely many
solutions in integers $x,y,z$ which satisfy the condition
$\op{GCD}(T_{x},T_{y})=1$?
\end{ques}

In the range $x<y<5\cdot 10^4$ there are exactly 39 solutions of our
equation and only one given by
\begin{equation*}
x=143,\quad y=237,\quad z=2301289,
\end{equation*} satisfied the
condition $\op{GCD}(T_{x},T_{y})=1$.

Unfortunately, we are unable to give an answer to the following
\begin{ques}
Does the equation $T_{x}^2+T_{y}^2=T_{z}^2$ have infinitely many
solutions in rational numbers?
\end{ques}

\section{Triangular numbers and palindromic numbers}\label{Section3}
Let us state the following
\begin{ques}
Let us fix $b\in\N_{>1}$. Is the set of triangular numbers which are
palindromic in base $b$ infinite?
\end{ques}

Let us remind that we say that the number $k$ is palindromic in base
$b$ if in the system with base $b$
\begin{equation*}
k=\sum_{i=0}^{m}a_{i}b^{i}
\end{equation*}
we have $a_{i}=a_{m-i}$ for $i=0,1,\ldots,\;m.$

In connection with this question we can prove the following
\begin{thm}
If $b=2,3,5,7,9$, then there are infinitely many triangular numbers
which are palindromic in base $b$.
\end{thm}
\begin{proof}
If $b=2$, then for $n=2^{2^k}+1$ we have
\begin{equation*}
t_{n}=(2^{2^{k}}+1)(2^{2^k-1}+1)=2^{2^{k+1}-1}+2^{2^{k}}+2^{2^{k-1}}+1=1\underbrace{00\ldots
00}_{k-zeros}11\underbrace{00\ldots 00}_{k-zeros}1_{2},
\end{equation*}
which proves that the number $t_{k}$ is palindromic in base 2.

If $b=3$ then we define $n=(3^k-1)/2$ and we get number
\begin{equation*}
t_{n}=\frac{3^{2k}-1}{3^2-1}=3^{2k-2}+3^{2k-4}+\ldots 3^{2}+3^{0},
\end{equation*}
which is clearly palindromic in base 3. Let us note that the number
$t_{n}$ is palindromic in base $b=9$ due to the identity
\begin{equation*}
t_{n}=\frac{9^{k}-1}{9-1}=11\ldots 11_{9}.
\end{equation*}

In the case when $b=5$ we put $n=(5^k-1)/2$ and we get the number
\begin{equation*}
t_{n}=\frac{5^{2k}-1}{8}=3\cdot5^{2k-2}+3\cdot5^{2k-4}+\ldots+3\cdot5^{2}+3\cdot5^{0},
\end{equation*}
which is palindromic in base 5.

If now $b=7$ then we define $n=(7^k-1)/2$ and we get the number
\begin{equation*}
t_{n}=\frac{7^{2k}-1}{8}=6\cdot7^{2k-2}+6\cdot7^{2k-4}+\ldots+6\cdot7^{2}+6\cdot7^{0},
\end{equation*}
which is palindromic number in base 7.
\end{proof}

\begin{rem}
{\rm Unfortunately, we are unable to prove that the set of
palindromic triangular numbers in base 10 is infinite. Let us note
that in the range $n<10^6$ there are exactly 35 values of $n$ with
such a property that the number $t_{n}$ is palindromic. However it
is easy to see that there are infinitely many triangular numbers
which are "almost" palindromic. This mean that at least one of the
numbers  $t_{n}\pm 1$ is palindromic. More precisely, if $n=2\cdot
10^{k+1}+1$ then
\begin{equation*}
t_{n}+1=2\underbrace{00\ldots 00}_{k-zeros}3\underbrace{00\ldots
00}_{k-zeros}2.
\end{equation*}
If we take $n=2\cdot 10^k+2$ then we have
\begin{equation*}
t_{n}-1=2\underbrace{00\ldots 00}_{k-1-zeros}5\underbrace{00\ldots
00}_{k-1-zeros}2.
\end{equation*}

}
\end{rem}

\section{Arithmetic progressions}\label{Section4}

In this section we consider the problem of construction of integer
solutions of certain diophantine equations connected with values of
some functions involving triangular and tetrahedral numbers in
arithmetic progressions.

We start with the following
\begin{thm}
The equation
\begin{equation*}
\frac{1}{t_{x}}+\frac{1}{t_{y}}=\frac{2}{t_{z}}
\end{equation*}
has infinitely many non-trivial solutions in positive integers
$x,y,z$. In other words: there are infinitely many three term
arithmetic progressions consisted of the numbers
$1/t_{x},1/t_{z},1/t_{y}$.
\end{thm}
\begin{proof}
In order to prove our theorem let us put $z=(y-x-1)/2$. Then we have
an equality
\begin{equation*}
t_{(y-x-1)/2}-\frac{2t_{x}t_{y}}{t_{x}+t_{y}}=\frac{(x+y+1)^2f(x,y)}{8(x^2+y^2+x+y)},
\end{equation*}
where $f(x,y)=-x+x^2-y-4xy+y^2$. Note that $f(1,5)=0$ and that the
following identity holds:
\begin{equation*}
f(x,y)=f(y,4y-x+1).
\end{equation*}
Form the above we can deduce that if we define:
$x_{0}=1,\;x_{1}=5,\;x_{n}=4x_{n-2}-x_{n-2}+1$, then for each $n$ we
have $f(x_{n},x_{n+1})=0$ and additionally, if $n\equiv 1\pmod{2}$,
then the number $(x_{n+1}-x_{n}-1)/2$ is integer. This conclusion
finishes the proof of our theorem.

In particular we have
\begin{equation*}
\frac{1}{t_{76}}+\frac{1}{t_{285}}=\frac{2}{t_{104}},\quad
\frac{1}{t_{1065}}+\frac{1}{t_{3976}}=\frac{2}{t_{1455}},\quad
\frac{1}{t_{14840}}+\frac{1}{t_{55385}}=\frac{2}{t_{20272}}
\;\ldots\;.
\end{equation*}

\end{proof}

\begin{thm}
The diophantine equation $z^2=(T_{x}+T_{y})/2$ has infinitely many
non-trivial solutions in integers $x,y,z$. In other words: there are
infinitely many three term arithmetic progressions consisted of the
numbers $T_{x},z^2,T_{y}$.
\end{thm}
\begin{proof}
Proof of our theorem is an immediate consequence of the identity
\begin{equation*}
\frac{T_{(u^2-1)/3}+T_{(2u^2-5)/3}}{2}=T_{u-1}^2
\end{equation*}
and the fact that for $u\equiv 1,2\pmod{3}$ the values of the
polynomials $(u^2-1)/3,\;(2u^2-5)/3$  are integers.

Let us note that we proved something more. Indeed, we have proved
that there are infinitely many three term arithmetic progressions
consisted of the numbers $T_{x},T_{u}^2,T_{y}$.
\end{proof}

\begin{thm}
The diophantine equation $t_{z}=\frac{x^4+y^4}{2}$ has infinitely
many non-trivial solutions in integers $x,y,z$. In other words:
there are infinitely many three term arithmetic progressions
consisted of the numbers $x^4,t_{z},y^4$.
\end{thm}
\begin{proof}
By a non-trivial solution of the equation $t_{z}=\frac{x^4+y^4}{2}$
we mean the solution $x,\;y,\;z$ which is not of the form
$x=m,y=m^2,z=m^4$ for certain $m\in\N$.

In order to prove our theorem it is enough to show that the
diophantine equation $(\star)\; u^2=4x^4+4y^4+1$ has infinitely many
solutions in integers. This is an easy consequence of the fact that
each triple of integers $(x,y,u)$ which is a solution of the
equation $(\star)$ give us a triple of integers $(x,y,(u-1)/2)$
which is a solution of the equation $t_{z}=(x^4+y^4)/2$.

Let us note the following identity
\begin{equation*}
(130w^2-128w+33)^2=4(60w^2-61w+16)^2+4(5w-2)^4+1.
\end{equation*}
This identity shows that the set of three terms arithmetic
progression consisted of the numbers $x^2,t_{z},y^4$ is infinite.
Thus, we can see that our theorem will be proved if we were able to
prove that the diophantine equation $v^2=60w^2-61w+16=:f(w)$ has
infinitely many solutions in integers. Let us note that $f(0)=4^2,$
and next
\begin{equation*}
(1921v + 14880w -7564)^2-f(248v + 1921w-976)=v^2-f(w).
\end{equation*}
Form the above identity we can deduce that if we define the
sequences $v_{n},w_{n}$ recursively by the equations
\begin{equation*}
\begin{cases}
w_{0}=0,\quad v_{0}=4,\\
w_{n}=248v_{n-1} + 1921w_{n-1}-976,\\
v_{n}=1921v_{n-1} + 14880w_{n-1}-7564,
\end{cases}
\end{equation*}
then for each $n\in\N$ we have the identity $v_{n}^2=f(w_{n})$. This
means that the equation $t_{z}=(x^4+y^4)/2$ has infinitely many
solutions $x_{n},y_{n},z_{n}$ given by
\begin{equation*}
x_{n}=v_{n},\quad y_{n}=60w_{n}^2-61w_{n}+16,\quad
z_{n}=65w_{n}^2-64w_{n}+16,\quad\quad n=0,1,2,\ldots\;.
\end{equation*}
In particular we have
$t_{16}=(4^4+2^4)/2,\;t_{15632}=(120^4+78^4)/2,\ldots \;$.
\end{proof}

In the light of the above theorem it is natural to state the
following
\begin{ques}
Does the equation  $t_{z}=x^4+y^4$ have infinitely many solutions in
integers?
\end{ques}
In the range $x<y<10^5$ there are two solutions of this equation.
There are the following triples:
$x=15,y=28,z=1153$;$\;x=3300,y=7712,z=85508608$.

\section{Varietes}\label{Section5}

We start with the following
\begin{thm}\label{iloraztetraedralnych}
There are infinitely many triangular numbers which are quotients of
tetrahedral numbers.
\end{thm}
\begin{proof}
It is easy to check that if
\begin{equation*}
x(u)=u,\quad y(u)=\frac{u^3+u^2+2u-4}{2},\quad z(u)=y(u)
\end{equation*}
or
\begin{equation*}
 x(u)=u,\quad y(u)=3T_{u},\quad
 z(u)=\frac{u^3+u^2+2u+2}{2},
\end{equation*}
then we have
\begin{equation*}
t_{z(u)}=\frac{T_{y(u)}}{T_{x(u)}}.
\end{equation*}
This ends proof of our theorem.
\end{proof}

In the light of the above theorem it is natural to ask about
triangular numbers which are products of tetrahedral numbers. We can
prove the following

\begin{thm}\label{iloczyntetraedralnych}
There are infinitely many triangular numbers which are product of
two tetrahedral numbers.
\end{thm}
\begin{proof}
We are interested in the integers $x,y,z$ which satisfy the equation
$t_{z}=T_{x}T_{y}$. In the table below we can find integer valued
polynomials $x_{i},y_{i},z_{i}$ which satisfy the equation
$t_{z_{i}}=T_{x_{i}}T_{y_{i}}$ for $i=1,\ldots,9$.
\begin{equation*}
\begin{array}{|c|l|}
\hline
  x_{1}(u) & 9u\\
  y_{1}(u) & (81u^3+27u^2+2u-2)/2=:f(u) \\
  z_{1}(u) & (f(u)+2)f(u)             \\
  \hline
  x_{2}(u) & 9u\\
  y_{2}(u) & 4f(u)+1\\
  z_{2}(u) & u(9u+1)(9u+2)(162u^3+54u^2+4u-3)\\
  \hline
  x_{3}(u) & 9u\\
  y_{3}(u) & 4f(u)+5\\
  z_{3}(u) & u(9u+1)(9u+2)(162u^3+54u^2+4u+3)\\
  \hline
  \hline
  x_{4}(u) & 9u-1\\
  y_{4}(u) & (81u^3-u-2)/2=:g(u)\\
  z_{4}(u) & g(u)(g(u)+2)\\
 \hline
  x_{5}(u) & 9u-1\\
  y_{5}(u) & 4g(u)+1\\
  z_{5}(u) & u(9u-1)(9u+1)(162u^3-2u-3)\\
 \hline
  x_{6}(u) & 9u-1\\
  y_{6}(u) & 4g(u)+5\\
  z_{6}(u) & u(9u-1)(9u+1)(162u^3-2u+3)\\
 \hline
 \hline
  x_{7}(u) & 9u-2\\
  y_{7}(u) & (81u^3-27u^2+2u-2)/2=:h(u)\\
  z_{7}(u) & h(u)(h(u)+2)\\
 \hline
  x_{8}(u) & 9u-2\\
  y_{8}(u) & 4h(u)+1\\
  z_{8}(u) & u(9u-2)(9u-1)(162u^3-54u^2+4u-3)\\
 \hline
  x_{9}(u) & 9u-2\\
  y_{9}(u) & 4h(u)+5\\
  z_{9}(u) & u(9u-2)(9u-1)(162u^3-54u^2+4u+3)\\
 \hline
\end{array}
\end{equation*}

\end{proof}

\begin{thm}
The diophantine equation $t_{p}^2+t_{q}^2=t_{r}^2+t_{s}^2$ has
infinitely many non-trivial solutions in integers.
\end{thm}
\begin{proof}
In order to prove our theorem it is enough to show that the
diophantine equation $(x^2-1)^2+(y^2-1)^2=(u^2-1)^2+(v^2-1)^2$ has
infinitely many solutions in odd integers. We use the method which
is very close to the method employed by Euler during his
investigation of integer solutions of the diophantine equation
$p^4+q^4=r^4+s^4,$ \cite[page 90]{Mor}.

Let us define  $f(x,y)=(x^2-1)^2+(y^2-1)^2$ and note that if
\begin{equation*}
x=T+c,\quad y=bT-d,\quad u=T+d,\quad v=bT+c,
\end{equation*}
then $f(x,y)-f(u,v)=-2a_{1}T-6a_{2}T^2-4a_{3}T^3=:g(T)$, where
\begin{align*}
&a_{1}=(b-1)c(c^2-1)-(b+1)d(d^2-1),\\
&a_{2}=(b^2-1)(c^2-d^2),\\
&a_{3}=c(b^3-1)+d(b^3+1).
\end{align*}
If we put $c=-b^3-1,\;d=b^3-1$, then $a_{3}=0$ and the equation
$g(T)=-8b^3(b^2-1)T(3T+b^4-2b^2-2)=0$ has two rational roots: $T=0$
and
\begin{equation*}
T=-(b^4-2b^2-2)/3.
\end{equation*}
Using the values of $T$ we have found and remembered that
$t_{-n}=t_{n-1}$ we can find solutions of the equation from the
statement of our theorem in the form
\begin{center}
\fboxrule=0mm
\begin{tabular}{ll}
\fbox{$p(b)=\frac{(b+1)(b^3+4b^2+2b+2)}{6}$}, & \fbox{$q(b)=\frac{b^5+b^3-2b-6}{6}$}\\
\fbox{$r(b)=\frac{(b+1)(b^3-4b^2+2b-2)}{6}$}, & \fbox{$s(b)=\frac{b(b^2-1)(b^2+2)}{6}$}.\\
\end{tabular}
\end{center}

It is easy to see that if $b\equiv 1\pmod{3}$, then the numbers
$p(b),q(b),r(b),s(b)$ are integers.
\end{proof}

W. Sierpi\'{n}ski in the paper \cite{Sier4} proved that the equation
$z^2=T_{x}+T_{y}$ has infinitely many solutions in integers. Next
step is the question if similar result can be proved if a cube
instead of a square is considered. As we will see the answer on such
modified question is affirmative

\begin{thm}
There are infinitely many cubes which are sums of two tetrahedral
numbers.
\end{thm}
\begin{proof}
Let us note the following equality
\begin{equation*}
\Big(\frac{x+6y}{2}\Big)^3-T_{x+5y-1}-T_{y-1}=\frac{1}{24}(x+6y)F(x,y),
\end{equation*}
where $F(x,y)=x^2-24y^2-4$. In order to finish the proof we must
show that the diophantine equation $F(x,y)=0$ has infinitely many
solutions in positive integers. In order to prove this let us note
that $F(2,0)=0$ and next that
\begin{equation*}
F(5x+24y,x+5y)=F(x,y).
\end{equation*}
Thus we can see that if we define
\begin{equation*}
x_{0}=2,\quad y_{0}=0,\quad x_{n}=5x_{n-1}+24y_{n-1},\quad
y_{n}=x_{n-1}+5y_{n-1},
\end{equation*}
then for each $n\in\N$ we have $2|x_{n}$ and $F(x_{n},y_{n})=0$.
This shows that the equation $z^3=T_{x}+T_{y}$ has infinitely many
solutions in integers. In particular we have:
$T_{1}+T_{19}=11^3,\;T_{19}+T_{197}=109^3,\;T_{197}+T_{1959}=1079^3\;\ldots\;
.$
\end{proof}

In the light of the result of Sierpi\'{n}ski and the above theorem
it is natural to ask the following
\begin{ques}
For which $n$ the diophantine equation $z^n=T_{x}+T_{y}$ has a
solution in positive integers?
\end{ques}
If $n=4$ then the smallest solution of this equation is:
$x=8,\;y=38,\;z=10$. Let us note that in the range $x<y<10^4$ this
equation has exactly six integer solutions.

\begin{thm}
There are infinitely many pairs of different tetrahedral numbers
with such a property that their product is a square of integer.
\end{thm}
\begin{proof}
We consider the problem of the existence of integer solutions of the
diophantine equation $z^2=T_{x}T_{y}$. In order to prove our theorem
it is enough to show that the equation $v^2=(x+2)(2x+1)/9=f_{1}(x)$
has infinitely many solutions in integers. Indeed, this is simple
consequence of the identity $T_{x}T_{2x}=x^2(x+1)^2f_{1}(x)$. We can
easily check the identity
\begin{equation*}
(8u + 17v + 10)^2-f_{1}(17u + 36v + 20)=v^2-f_{1}(u).
\end{equation*}
Because $f_{1}(1)=1$ we can see that if we define sequences
$u_{n},v_{n}$ recursively in the following way
\begin{equation*}
u_{0}=1,\quad v_{0}=1,\quad u_{n}=17u_{n-1} + 36v_{n-1} + 20,\quad
v_{n}=8u_{n-1} + 17v_{n-1} + 10,
\end{equation*}
then for each  $n\in\N$ the following identity holds
\begin{equation*}
(v_{n}u_{n}(u_{n}+1))^2=T_{u_{n}}T_{2u_{n}}.
\end{equation*}
In particular $189070^2=T_{73}T_{146},\;
7559616818^2=T_{2521}T_{5042},\;\ldots\;. $

At the end let us note that the proof of our theorem can be
performed with the use of the identity
$T_{x}T_{2x+2}=(x+1)^2(x+2)^2f_{2}(x)$, where $f_{2}(x)=x(2x+3)/9$.
It is easy to prove that the equation $v^2=f_{2}(u)$has infinitely
many solutions in integers.

In the light of the proof of our theorem and the remark above we can
state an interesting question concerning the existence of infinite
set of triples of integers $x,y,z$ which satisfy the equation
$z^2=T_{x}T_{y}$ with the condition $2x+2<y$?
\end{proof}
\begin{rem}
{\rm We do not know, if there exist three different tetrahedral
numbers in geometric progression, but it is possible to prove that
there are infinitely non-trivial rational solutions of the
diophantine equation $T_{x}T_{y}=T_{z}^2$. Proof of this fact can be
found in \cite{Ulas}.

}
\end{rem}

\bigskip

 \hskip 4.5cm       Maciej Ulas

 \hskip 4.5cm       Jagiellonian University

 \hskip 4.5cm       Institute of Mathematics

 \hskip 4.5cm       {\L}ojasiewicza 6

 \hskip 4.5cm       30 - 348 Krak\'{o}w, Poland

 \hskip 4.5cm      e-mail:\;{\tt Maciej.Ulas@im.uj.edu.pl}

 \end{document}